\documentclass[12pt]{article} 

\usepackage{times}
\usepackage{latexsym}
\usepackage{amsmath}
\usepackage{amsfonts}
\usepackage{xcolor}
\usepackage{tikz}
\usepackage{amsthm}
\begingroup
    \makeatletter
    \@for\theoremstyle:=definition,remark,plain\do{%
        \expandafter\g@addto@macro\csname th@\theoremstyle\endcsname{%
            \addtolength\thm@preskip\parskip
            }%
        }
\endgroup
\makeatletter
\newtheorem*{rep@theorem}{\rep@title}
\newcommand{\newreptheorem}[2]{%
\newenvironment{rep#1}[1]{%
 \def\rep@title{#2 \ref{##1}}%
 \begin{rep@theorem}}%
 {\end{rep@theorem}}}
\makeatother
\usepackage{parskip}
\usepackage{GWmacros}
\usepackage{epstopdf}

\newreptheorem{Theorem}{Theorem}

\title{The weak Bruhat order for random walks on Coxeter groups}
\author{Graham White}
\date{\today}
\begin{document}

\maketitle


\begin{abstract}
We show that for the simple random walk on a Coxeter group generated by the Coxeter generators and identity, the likelihoods of being at any pair of states respect the weak Bruhat order. That is, after any number of steps, the most likely element is the identity, probabilities decrease along any geodesic from the identity, and the least likely element is the longest element, if the group is finite. The result remains true when different generators have different probabilities, so long as the identity is at least as likely as any other.
\end{abstract}

\section{Introduction}

If a random walk is started in a known state and run for several steps, one may examine the probabilities that it is now in each possible state. A likelihood order is a partial order on the state space, so that if one state is larger than another, then the random walk is always more likely to be in the former state than the latter, after any number of steps. The main result of this paper is that for any Coxeter system $(W,S)$, the weak Bruhat order is a likelihood order for the simple random walk on $W$ generated by $S \cup \{1_W\}$. 

\begin{Theorem}
\label{the:main}
For any Coxeter system $(W,S)$, consider the simple random walk on $W$, starting at the identity, and at each step multiplying on the right by an element of $S$ or by the identity, each with probability $\frac{1}{|S|+1}$. Then for any $n$, and any two states $w$ and $w'$, if $w \leq_B w'$ then the probability that the random walk is at $w$ after $n$ steps is at least the probability that it is at $w'$. 
\end{Theorem}

In type A, this result describes a partial order for the (appropriately lazy) adjacent transposition walk on the symmetric group. For analysis of the mixing time of this walk, see Section 4 of \cite{Compgroups}. The adjacent transposition walk is a special case of the interchange process \cite{Interchange}.

Likelihood orders can describe the most and least likely states of a random walk. In particular, if a random walk has uniform stationary distribution, then the separation distance from the stationary distribution depends only on the probability of being at the least likely state. Thus, knowing which state is the least likely, together with a lower bound on the probability of being at that state, produces an upper bound on the separation distance mixing time. Upper bounds on separation distance give upper bounds on total variation distance, so upper bounds on total variation mixing times follow. 

Lower bounds on mixing times are obtained by analysing a set of unlikely states. Knowledge of which states are the least likely via a likelihood order can inform the choice of such a set. For instance, see Section 4 of \cite{Thumb}.

Some results regarding likelihood orders for various random walks on the symmetric group $S_n$ are given in \cite{MeganLikelihood} and \cite{MeganInvolutions}. These likelihood orders are shown to hold after enough steps (for example, after $O(n^2)$ steps). Likelihood orders have also been considered by Diaconis and Isaacs in \cite{DiaconisIsaacs}, where they prove that for any symmetric random walk on a group, after any even number of steps the most likely state is the initial one. They also give likelihood orders for several random walks on the cycle.

\section{Preliminaries}

In this section, some necessary background results are recalled. Section \ref{sec:cayley} defines Cayley graphs and a useful family of their symmetries, and Section \ref{sec:distances} discusses the sets of vertices in the Cayley graph which are closer to one end of a given edge than to the other. Section \ref{sec:bruhat} defines the weak Bruhat order.

A Coxeter system $(W,S)$ is a group $W$ together with a presentation of a certain form.

\begin{Definition}
A \emph{Coxeter presentation} is a presentation of the form $$\pres{s_1,s_2,\dots,s_n}{\{s_i^2\}_{i=1}^n, \{(s_is_j)^{m_{ij}}\}_{i \neq j}}$$ where each $m_{ij}, i \neq j$ is either a positive integer at least two, or $\infty$, indicating the lack of that relation. 
\end{Definition}

A good example of a Coxeter group is the symmetric group $S_n$, which has the Coxeter presentation $$S_n = \left\langle s_1, s_2, \dots, s_{n-1}\left|\begin{array}{lcl}
s_i^2 & \text{for each} & 1 \leq i \leq n-1 \\
(s_is_{i+1})^3 & \text{for each} & 1 \leq i \leq n-2 \\
s_is_js_i^{-1}s_j^{-1} & \text{if} & |i-j| > 1 \\
\end{array}\right.\right\rangle.$$ 

In this presentation, the generator $s_i$ represents the transposition $(i \; i+1)$. This presentation has $m_{ij} = 3$ when $|i-j|=1$ and $m_{ij} = 2$ for $|i-j|>1$.

\subsection{Cayley graphs}
\label{sec:cayley}

Given a group $W$ and a generating set $S$, the Cayley graph $\Gamma(W,S)$ is defined as follows

\begin{Definition}
The graph $\Gamma(W,S)$ has a vertex for each element of $W$, and for each $w \in W$ and each $s \in S$, there is an edge from $w$ to $ws$. It will often be convenient to label the edge $(w,ws)$ by the generator $s$.  
\end{Definition} 

In the present setting, groups will always be generated by elements of order two, so Cayley graphs will be undirected.

It will be necessary to have the following results regarding certain symmetries of Cayley graphs.

\begin{Definition}
Consider a Cayley graph $\Gamma(W,S)$. For any $x \in W$, let $L_x$ be the left multiplication map on $\Gamma(W,S)$ which takes $w$ to $xw$, for each $w \in W$.
\end{Definition}

\begin{Proposition}
For any Cayley graph $\Gamma(W,S)$ and any $x \in W$, the map $L_x$ is an automorphism of $\Gamma(W,S)$. Further, $L_x$ preserves the edge labels of $\Gamma(W,S)$.
\end{Proposition}
\begin{proof}
The map $L_{x^{-1}}$ is the inverse of $L_x$, so $L_w$ is a bijection. To check that $L_x$ preserves edges of $\Gamma(W,S)$, observe that for each edge $(w,ws)$ of $\Gamma(W,S)$, the image under $L_x$, $(xw,xws)$, is also an edge of $\Gamma(W,S)$, and that these two edges have the same label.
\end{proof}

Random walks on the group $W$ can be understood via the Cayley graph. In particular, if a random walk is defined by at each step multiplying by an element of $S$, then consider the set of paths in $\Gamma(W,S)$ of length $n$ which start at the identity. The probability $P^n(w)$ that the walk is at $w$ after $n$ steps is equal to the proportion of these paths which end at $w$. Lazy walks can be considered by including $1_W$ in $S$.

\subsection{Distances in $\Gamma(W,S)$}
\label{sec:distances}

It will be important to understand relative distances in the Cayley graphs of Coxeter groups. Consider the Cayley graph with the usual graph metric --- that is, each edge has length $1$, and the distance $d(w,x)$ between two vertices $w$ and $x$ is the number of edges in the shortest path connecting them. As usual in the theory of Coxeter groups, $l(w)$ will denote the distance from the identity $d(1_W,w)$. Equivalently, $l(w)$ is the fewest number of generators which can be multiplied to produce $w$.

For this section, fix $w \in W$ and $s \in S$, with $l(w) < l(ws)$. That $l(w) < l(ws)$ is not used in this section, but is consistent with how these results will be used in Section \ref{sec:main}.

\begin{Definition}
\label{def:colours}
Let $\Gamma(W,S)$ be the Cayley graph of a Coxeter system $(W,S)$. For the fixed adjacent vertices $w$ and $ws$ of $\Gamma(W,S)$, colour each vertex of $\Gamma(W,S)$ white if it is closer to $w$ than to $ws$ and black if it is closer to $ws$ than to $w$. 
\end{Definition}

\begin{Proposition}
Each vertex of $\Gamma(W,S)$ is coloured white or black, but not both.
\end{Proposition}
\begin{proof}
Each vertex of $\Gamma(W,S)$ has at most one colour, because it cannot be both closer to $w$ than to $ws$ and the reverse. To show that each vertex is coloured, it is sufficient to show that no vertex can be equidistant from $w$ and $ws$.

The Coxeter relations of $(W,S)$ each have even length, so $\Gamma(W,S)$ is a bipartite graph, and hence the distances from any vertex to $w$ and $ws$ have opposite parities. Thus each vertex is coloured, completing the proof.
\end{proof}

\begin{Definition}
Continuing from Definition \ref{def:colours}, colour grey each edge which connects a white vertex to a black vertex. 
\end{Definition}

\begin{Lemma}
\label{lem:greydist}
If $(x,xt)$ is a grey edge, with $x$ white and $xt$ black, then $d(x,w) = d(xt,ws)$. (The generator $t$ may be equal to $s$, but need not be.)
\end{Lemma}
\begin{proof}
The two vertices $w$ and $ws$ are adjacent, as are the vertices $x$ and $xt$. The vertex $x$ is white, and $xt$ is black. Thus, the following relations between distances hold 
\begin{align*}
d(x,w) + 1 &= d(x,ws) \\ 
d(xt,ws) + 1 &= d(xt,w) \\
d(x,ws) &= d(xt,ws) \pm 1\\
d(xt,w) &= d(x,w) \pm 1\\
\end{align*} 
Adding these four equations, each of the two $\pm$ signs must be a plus. Thus $d(x,w) = d(xt,ws)$, as required.
\end{proof}

\begin{Lemma}
\label{lem:greypath}
Under the conditions of Lemma \ref{lem:greydist}, $w^{-1}x = sw^{-1}xt$.
\end{Lemma}
\begin{proof}
Let $\boldsymbol{\omega}$ be a reduced word for $w^{-1}x$. From Lemma \ref{lem:greydist}, $l(sw^{-1}xt) = l(w^{-1}x)$. Thus, $s\boldsymbol{\omega} t$ is a word of length two greater than the minimum length of any equivalent word. By the deletion condition (Section 1.7 of \cite{Humphreys}), there is a reduced word for $sw^{-1}xt$ which can be obtained by deleting two letters from $s\boldsymbol{\omega} t$. However, Lemma \ref{lem:greydist} also implies that the words $s\boldsymbol{\omega}$ and $\boldsymbol{\omega} t$ are reduced, so the two letters removed from the word $s\boldsymbol{\omega} t$ must be the initial $s$ and the final $t$. Therefore $s\boldsymbol{\omega} t$ and $\boldsymbol{\omega}$ are equivalent words, so $w^{-1}x = sw^{-1}xt$.
\end{proof}

\begin{Proposition}
\label{prop:greyflip}
The map $L_{wsw^{-1}}$ interchanges the endpoints of any grey edge.
\end{Proposition}
\begin{proof}
Let $(x,xt)$ be a grey edge, with $x$ white and $xt$ black. The image $L_{wsw^{-1}}(x)$ is
\begin{align*}
L_{wsw^{-1}}(x) &= wsw^{-1}x \\
& = wssw^{-1}xt \text{   (by Lemma \ref{lem:greypath})} \\
& = xt \\
\end{align*}
The map $L_{wsw^{-1}}$ is an involution, so $L_{wsw^{-1}}(xt) = x$.
\end{proof}

\begin{Corollary}
Any vertex of $\Gamma(W,S)$ is adjacent to at most one grey edge.
\end{Corollary}
\begin{proof}
The function $L_{wsw^{-1}}$ is well defined. If any vertex were adjacent to more than one grey edge, then Proposition \ref{prop:greyflip} implies that $L_{wsw^{-1}}$ is multivalued, a contradiction.
\end{proof}

The results in this section are not new --- they are standard facts about the geometry of $\Gamma(W,S)$, proven again here to draw out the key pieces. The set of grey edges is commonly referred to as a wall, and could be defined as the set of edges preserved by the reflection $L_{wsw^{-1}}$.

\subsection{The weak Bruhat order}
\label{sec:bruhat}

Let $(W,S)$ be a Coxeter system. The (right) weak Bruhat order on $(W,S)$ is defined as follows (Chapter 3 of \cite{Bjorner}).

\begin{Definition}
Let $w$ and $w'$ be elements of $W$. Then $w \leq_B w'$ if there is a reduced word for $w$ which can be multiplied on the right by elements of $S$ to produce a reduced word for $w'$.	
\end{Definition}

An equivalent formulation in terms of the Cayley graph of $(W,S)$ is

\begin{Definition}
Let $w$ and $w'$ be elements of $W$. Then $w \leq_B w'$ if there is a minimal length path in $\Gamma(W,S)$ from the identity element $1_W$ to $w'$ which passes through $w$.	
\end{Definition}

These two definitions are equivalent because edges in the Cayley graph $\Gamma(W,S)$ correspond to multiplication on the right by an element of $S$.

\section{Proof of main result}
\label{sec:main}

The main result of this paper is that the weak Bruhat order arises as a likelihood order for random walks on $(W,S)$.

\begin{repTheorem}{the:main}
For any Coxeter system $(W,S)$, consider the simple random walk on $W$, starting at the identity, and at each step multiplying on the right by an element of $S$ or by the identity, each with probability $\frac{1}{|S|+1}$. Then for any $n$, and any two states $w$ and $w'$, if $w \leq_B w'$ then the probability that the random walk is at $w$ after $n$ steps is at least the probability that it is at $w'$. 
\end{repTheorem}

To prove Theorem \ref{the:main}, it suffices to consider $w$ and $w'$ which are adjacent --- that is, when $w' = ws$ for some $s \in S$, with $l(w) < l(ws)$. If the result is true for adjacent vertices, then the general case follows by induction. Thus, the theorem reduces to the following proposition.

\begin{Proposition}
\label{prop:likelihood}
Let $w \in W$ and $s \in S$, with $l(ws) > l(w)$. Then for any $n$, $P^n(ws) < P^n(w)$.
\end{Proposition}
\begin{proof}
Let $\cP$ be the set of paths of length $t$ from the identity $1_W$ to $w$, with each step being either the identity or an element of $S$. Let $\cP'$ be the set of paths of length $t$ from $1_W$ to $ws$. To prove this proposition, it suffices to construct an injection from $\cP'$ to $\cP$. Consider an element $\alpha$ of $\cP'$, and write it as $$\alpha = (1_W=a_0,a_1,a_2,\dots,a_n=ws).$$

The path $\alpha$ starts at a point closer to $w$ than to $ws$, because $l(ws) > l(w)$, and $\alpha$ ends at $ws$, which is closer to $ws$ than to $w$. Thus at some point, $\alpha$ crosses from a vertex closer to $w$ to a vertex closer to $ws$ --- that is, this path crosses a grey edge. Given $\alpha$, let $i$ be the last time at which $\alpha$ either just crossed a grey edge or just stayed in place on an endpoint of a grey edge. Using $\oplus$ to denote concatenation, define the sequence of vertices 
\begin{align*}
f(\alpha) &= (a_j)_{j=0}^{i-1} \oplus (L_{wsw^{-1}}(a_j))_{j=i}^{n} \\
 &= (1_W=a_0,a_1,\dots,a_{i-1},L_{wsw^{-1}}(a_i),L_{wsw^{-1}}(a_{i+1}),\dots,L_{wsw^{-1}}(a_n)) = w.\\
\end{align*}

That is, the sequence $f(\alpha)$ is identical to $\alpha$ up until time $i-1$, and from time $i$ onwards, it is reflected by $L_{wsw^{-1}}$. 

\begin{Proposition}
The sequence $f(\alpha)$ is a path. That is, each two consecutive entries are either adjacent in the graph $\Gamma(W,S)$, or equal.
\end{Proposition}
\begin{proof}
It must be checked that each two consecutive entries in $f(\alpha)$ are either adjacent in the graph $\Gamma(W,S)$, or equal. This is immediate for each consecutive pair except for $a_{i-1}$ and $L_{wsw^{-1}}(a_i)$.

By the definition of $i$, the vertex $a_i$ is a black vertex adjacent to exactly one grey edge, and $a_{i-1}$ is one of the two endpoints of that edge. From Proposition \ref{prop:greyflip}, $L_{wsw^{-1}}(a_i)$ is either equal to $a_{i-1}$ or connected to $a_{i-1}$ by this grey edge.
\end{proof}

The map $L_{wsw^{-1}}$ interchanges $ws$ and $w$, so $f$ is a function from $\cP'$ to $\cP$. All that remains is to show that $f$ is an injection. The function $f$ is an involution, because it applies $L_{wsw^{-1}}$ to the part of the path from time $i$ onwards, the map $L_{wsw^{-1}}$ is an involution, and moving from $\alpha$ to $f(\alpha)$ does not change the definition of $i$ (but rather interchanges the two cases in the definition of $i$).

Thus $f$ is an involution from $\cP'$ to $\cP$, so there are at least as many paths of length $n$ from $1_W$ to $w$ as from $1_W$ to $ws$, completing the proof of Proposition \ref{prop:likelihood}.
\end{proof}

\begin{Corollary}
For any Coxeter system $(W,S)$ and the corresponding random walk described by Theorem \ref{the:main}, the most likely element after $n$ steps is the identity. If $W$ is finite, then the least likely element is the longest element. Here, most and least likely are not necessarily strict.
\end{Corollary}

\begin{Example}
Consider the random walk on $S_4$ generated by the adjacent transpositions $(1 \; 2)$, $(2 \; 3)$, and $(3 \; 4)$, as well as the identity. For any $n$, the most likely element after $n$ steps is the identity and the least likely is the reversal $(1 \; 4)(2 \; 3)$. Theorem \ref{the:main} does not address the relative likelihoods of the transpositions $s = (1 \; 2)$ and $t = (2 \; 3)$, but both are more likely than $(1 \; 3) = sts = tst$.
\end{Example}

\begin{Example}
Take the simple random walk on the square grid $\Z \times \Z$ generated by $(\pm 1,0)$, $(0,\pm 1)$ and $(0,0)$. This is a relabelling of the Cayley graph of the Coxeter group $$D_\infty \times D_\infty = \pres{s,t,u,v}{s^2,t^2,u^2,v^2,susu,svsv,tutu,tvtv}.$$ For any $n$, the most likely element after $n$ steps is the identity $(0,0)$, and the next most likely are the four adjacent vertices (which are equally likely, by symmetry). The vertex $(1,3)$ is always more likely than $(2,3)$, but Theorem \ref{the:main} doesn't address the relative likelihoods of $(1,3)$ and $(2,2)$.
\end{Example}

\begin{Example}
Finally, consider the simple random walk on a $d$--regular tree, with laziness $\frac{1}{d+1}$. This is the Cayley graph of the free product of several copies of the group $\Z / 2\Z$. For any $n$, the most likely vertex after $n$ steps is the initial vertex, the next most likely are the adjacent vertices, then the vertices at distance two, and so on.
\end{Example}

Theorem \ref{the:main} may be strengthened to not require that all generators have equal probability.

\begin{Theorem}
For any Coxeter system $(W,S)$, consider a random walk on $W$, starting at the identity, and at each step multiplying on the right by an element $s \in S$ or by the identity, with probabilities $p_s$ or $p_{\text{id}}$. As long as each $p_s$ is less than $p_{\text{id}}$, the conclusion of Theorem \ref{the:main} holds.
\end{Theorem}
\begin{proof}
The only part of the proof that must be changed is the proof of Proposition \ref{prop:likelihood}, comparing the probabilities of the states $ws$ and $w$, for $w \in W$ and $s \in S$ with $l(ws) > l(w)$. With this choice of $s$, divide the probability of multiplying by $1_W$ into two parts, of probabilities $p_s$ and $p_{\text{id}} - p_s$. Where the proof of Proposition \ref{prop:likelihood} pairs up the events of multiplying by $s$ or by $1_W$, use only the first of these parts, which is an event of equal probability to that of multiplication by $s$.
\end{proof}

\bibliographystyle{hplain}
\bibliography{bib}
\end{document}